\documentclass{amsart}

\newtheorem{thm}{Theorem}[section]
\newtheorem{lem}[thm]{Lemma}

\newtheorem{DFN}[thm]{Definition}
\newtheorem{RMK}[thm]{Remark}
\newtheorem{EG}[thm]{Example}
\newenvironment{dfn}{\begin{DFN} \rm}{\end{DFN}}

\newenvironment{eg}{\begin{EG} \rm}{\end{EG}}

\def\dis{\displaystyle}
\def\cal{\mathcal}
\def\cf{{\it cf.\ }}
\def\ie{{\it i.e.\ }}
\def\C{{\mathbb C}}

\def\la{\lambda}
\def\m{\mu}

\def\z{\zeta}
\def\o{\omega}

\def\inv{\mbox{inv}}

\def\pf{\par \noindent {\em Proof.}\hspace{.5em}}
\def\endpf{ \framebox(5,7)[bl]{} \par \bigskip \par}

\begin{document}

\title{Canonical Resolution of a Quasi-ordinary Surface Singularity}
\author{Chunsheng Ban \and  Lee J. McEwan}

\begin{abstract}
We describe the embedded resolution of a quasi-ordinary surface singularity
$(V,p)$ which results from applying the canonical resolution of
Bierstone-Milman to $(V,p)$. We show that this process depends solely on
the characteristic pairs of $(V,p)$, as predicted by Lipman. We describe
the process explicitly enough that a resolution graph for $f$ could in
principle be obtained by computer using only the characteristic pairs.
\end{abstract}

\maketitle

\bigskip

\noindent {\bf INTRODUCTION.} The Jungian approach to resolving the
singularities of an embedded surface $V\subset \C^3$ begins with a
projection $\pi\!: V\longrightarrow\C^2$
with (reduced) discriminant locus $\Delta\subset\C^2$. Let
$\sigma\!:(M,\sigma^{-1}(\Delta))\longrightarrow(\C^2,\Delta)$
be an embedded resolution of the plane curve $\Delta$. Thus, $M$ is
smooth (of dimension $2$) and $\sigma^{-1}(\Delta)$ has only normal
crossings. The bimeromorphic map $\sigma$ and the projection
$\pi$ induce a bimeromorphic map
$\sigma'\!: V'=V\times_{\C^2}M\longrightarrow V$, and a projection
$\pi'\!:V'\longrightarrow M$ whose discriminant locus has only normal
crossings. We have the following diagram.

\vspace{.8in}
\hspace{1in}
\begin{picture}(60,25)(0,0)
\setlength{\unitlength}{1mm}
\put(5,0){$M$}
\put(12,1){\vector(1,0){30}}
\put(45,0){$\C^2$}
\put(5,20){$V'$}
\put(12,21){\vector(1,0){30}}
\put(45,20){$V$}
\put(3,13){$\pi'$}
\put(7,18){\vector(0,-1){13}}
\put(46,18){\vector(0,-1){13}}
\put(47,13){$\pi$}
\put(27,2){$\sigma$}
\put(27,22){$\sigma'$}
\put(13,7){$\sigma^{-1}(\Delta)$}
\put(56,7){$\Delta$}
\put(27,8){\vector(1,0){27}}
\put(17,6){\vector(-2,-1){6}}
\put(57,6){\vector(-2,-1){6}}
\end{picture}

\bigskip

\noindent This leads us to consider surface singularities in $\C^3$
admitting a local finite projection to $\C^2$ with normal crossings.
Such singularities are called quasi-ordinary.\footnote{Higher dimensional
quasi-ordinary singularities are similarly defined. In this paper we
only consider surface singularities in $\C^3$.} Lipman and others have
studied this kind of singularity in detail (\cf \cite{L1}, \cite{L2},
\cite{L3}, and \cite{G}). Each quasi-ordinary singularity is equipped
with a set of pairs $\{(\la,\m)\}$ of non-negative rational numbers,
called characteristic pairs. These pairs determine much of the geometry
and topology of the singularity.

One approach to understanding a quasi-ordinary singularity is to study
its resolution. In \cite{L1}, Lipman described a definite procedure for
resolving a quasi-ordinary surface singularity $(V,p)$ embedded in $\C^3$
(for an outline of this procedure, see \cite{L2}, pages 164--165 and pages
169--171). He shows that at each stage of the resolution, the singular
locus remains quasi-ordinary, and the characteristic pairs at each stage
are determined by those of the previous stage and the process employed.
In this paper we study the {\em embedded resolution\/} of a quasi-ordinary
surface singularity. Recently, several constructive algorithms
for canonical embedded resolution of singularities have emerged; see
Bierstone and Milman (\cite{BM1}, \cite{BM2}), Villamayor (\cite{EV}),
and Moh (\cite{M}). In this note we apply Bierstone and Milman's
algorithm. Specifically, we replace Lipman's procedure by applying Bierstone
and Milman's algorithm for canonical embedded resolution, and show that
at each stage of the algorithm, the center to blow up is
determined by the characteristic pairs and the process employed. Thus we establish a canonical embedded resolution for quasi-ordinary surfaces; the
fact that the data of the resolution depend only on characteristic pairs
will be used in a sequel to prove the existence of simultaneous embedded
resolution for equisingular families.

The authors are thankful for several conversations and/or correspondences
with E.\ Bierstone, G.\ Kennedy, J.\ Lipman, O.\ Villamayor and S.\ Encinas.
A partial version of these results were obtained earlier in terms of the
algorithm of \cite{V}, \cite{EV}. We ultimately found it easier to use the
more explicit bottom-up approach of Bierstone and Milman's treatment of
canonical resolution, but we learned this subject equally from E.\ Bierstone,
O.\ Villamayor and S.\ Encinas. We gratefully acknowledge our debt
to them. We are told that S.\ S.\ Abhyankar has presented a proof of our main
result in a new appendix to the second edition of his book \cite{A}, which
we have not seen, but which may not contain an explicit recipe in terms of
the characteristic pairs.

\section{Preliminaries} \label{pre}

\subsection{Quasi-ordinary Singularities}
A quasi-ordinary singularity is an analytic germ $(V,p)$ of
dimension $d$ which admits a finite map of analytic germs
$\pi\!: (V,p) \to (\C^d,0)$ whose discriminant locus $\Delta$
(the hypersurface in $\C^d$ over which $\pi$ ramifies) has only
normal crossings. We consider a quasi-ordinary surface singularity
embedded in $(\C^3,0)$ here. In this case, we can choose local
coordinates $x$, $y$, $z$ such that $\pi(x,y,z)=(x,y)$ and such
that $(V,p)$ is defined by an equation
\[ f(x,y,z)=z^m+c_{m-1}(x,y)z^{m-1}+\dotsb +c_0(x,y)=0 \]
where the $c_i$ are power series. Then $(V,p)$ being
quasi-ordinary means that the {\em discriminant\/} $D(x,y)$
of $f$ (considered as a polynomial in $z$) has the form
\[ D(x,y)=x^ay^b\epsilon(x,y), \quad \epsilon (0,0)\neq 0. \]
It is known that the roots of $f$ are fractional power series.
Let $\z=H(x^{1/n},y^{1/n})$ be one root of $f$ ($H(x,y)$ a power
series). The other roots (the conjugates of $\z$) can be expressed as
\[ \z_i=H_i(x^{1/n},y^{1/n})=H(\o_{i1}x^{1/n},\o_{i2}y^{1/n}) \]
with the $\o_{ij}$ the $n$-th roots of unity. By letting $\z_1=\z$,
the polynomial $f$ can be written as
\begin{eqnarray}
f(x,y,z) & = & z^m+c_{m-1}(x,y)z^{m-1}+\dotsb +c_0(x,y) \nonumber \\
         & = & \prod_{i=1}^{m}(z-\z_i),
\end{eqnarray}
so that the $c_i$ are the elementary symmetric functions of
$\z_i$. We call $\z$ a parametrization of $(V,p)$.

Since $x^ay^b\epsilon(x,y)=D(x,y)=\prod_{i\neq j}(\z_i-\z_j)$,
we have
\[ \z_i-\z_j=x^{u/n}y^{v/n}\epsilon_{ij}(x^{1/n},y^{1/n}), \quad
  \epsilon_{ij}(0,0)\neq 0  \]
for some $u$ and $v$ (depending on $i$, $j$). The fractional
monomials $M_{ij}=x^{u/n}y^{v/n}$ so obtained are called the
{\em characteristic monomials\/} of $\z$, and the exponents
$(u/n,v/n)$ are called the {\em characteristic pairs\/} of $\z$. These
pairs satisfy certain conditions (\cf \cite{L2}, Proposition 1.5); for
example, they are totally ordered by $(\la_1,\m_1)\leq(\la_2,\m_2)$
if and only if $x^{\la_1}y^{\m_1}$ divides $x^{\la_2}y^{\m_2}$ (\ie
$\la_1\leq\la_2$ and $\m_1\leq\m_2$). It turns out that these
pairs determine quite a lot of the geometry and topology of
$(V,p)$ (\cf \cite{L1}, \cite{L2}, \cite{L3}, and \cite{G}).

Although a quasi-ordinary singularity may have different characteristic
pairs which depend on the choice of $\z$, the characteristic pairs
of a {\em normalized\/} parametrization determine and are determined
by the local topological type of $(V,p)$ (\cf \cite{G}); that is,
there is a set of characteristic pairs which are independent of the
choice of $\z$.

\begin{dfn}
A parametrization
\[ \z=x^{a/n}y^{b/n}H(x^{1/n},y^{1/n}), \]
where $H(0,0)\neq 0$ is {\em normalized\/} if
\begin{enumerate}
\item $a$ and $b$ are not both divisible by $n$,
\item if $a+b<n$, then both $a>0$ and $b>0$,
\item labeling the characteristic pairs
$(\la_i,\m_i)_{1\leq i\leq s}$ of $\z$ so that $\la_1\leq\la_2\leq
\cdots\leq\la_s$ and $\m_1\leq\m_2\leq\cdots\leq\m_s$, we have
$(\la_1,\la_2,\dots,\la_s)\geq (\m_1,\m_2,\dots,\m_s)$
(lexicographically).
\end{enumerate}
\label{normal}
\end{dfn}

\begin{lem}
Any parametrization $\z$ of a quasi-ordinary singularity $(V,p)$ can be
changed to a normalized one by the following type of change of coordinates:
\[ \begin{array}{crcl}
  & x' & = & x \\
(i) \quad & y' & = & y \\
  & z' & = & z-p(x,y)
   \end{array} \qquad
   \begin{array}{crcl}
  & x' & = & z \\
(ii) \quad & y' & = & y \\
  & z' & = & x
   \end{array} \qquad
   \begin{array}{crcl}
  & x' & = & y \\
(iii) \quad & y' & = & x \\
  & z' & = & z
   \end{array}.      \]
\label{trans}
   \end{lem}
\pf Let $p(x,y)$ be the sum of the terms $x^{\alpha}y^{\beta}$ in $\z$
with $\alpha$, $\beta$ both integers. By a change of type (i), we get a
parametrization $\z'=\z-p(x,y)$. By the property of characteristic pairs
(\cf \cite{L2}, Proposition 1.5), the smallest characteristic monomial
divides every term in $\z'$. So,
$\z'=x^{\la_1}y^{\m_1}H(x^{1/n},y^{1/n})$
for some $H$. Therefore, we may assume that
$\z=x^{\la_1}y^{\m_1}H(x^{1/n},y^{1/n})$, and 1 in (\ref{normal})
is satisfied. If 2 in (\ref{normal}) is not satisfied, from the proof
of the {\em inversion formula\/} (\cf \cite{G}, Appendix) a change of
type (ii) can take care of this. Upon a change of type (iii), 3 in
(\ref{normal}) can be obtained. \hfill\endpf

\subsection{Resolution of Quasi-ordinary Singularities}
In \cite{L1} J. Lipman described a definite procedure to resolve a
quasi-ordinary singularity (for an outline of this procedure, see
\cite{L2}, pages 164--165 and pages 169--171). This procedure determines
a tree of {\em special transforms\/} of a quasi-ordinary singularity
$(V,p)$. Each special transform in this tree is again a quasi-ordinary
singularity. The center to blow up at each stage is determined by the
characteristic pairs, and the characteristic pairs of the resulting
special transform are determined by the pairs in the previous stage
and the process employed (\cf \cite{L2}, Theorem~3.3). Let
\begin{equation}
\z=x^{\la}y^{\mu}H(x^{1/n},y^{1/n}), \quad H(0,0)\neq 0
\label{zeta}
\end{equation}
be a normalized parametrization of a quasi-ordinary singularity
$(V,p)$, where $p=(0,0,0)$ is the origin.
Roughly speaking, we blow up a {\em permissible curve\/} (an
equi-multiple curve) whenever it is possible, and a point
otherwise. The possible permissible curve is the plane section
cut out by coordinate plane $x=0$ (if $\la\geq 1$) or $y=0$
(if $\mu\geq 1$). By Bierstone-Milman's algorithm, however, we will at
sometime blow up a point even though there is a permissible curve.

The following lemma by Lipman shows how the characteristic pairs
change after blowing up the origin (quadratic transform) or a
permissible curve (monoidal transform). (\cf \cite{L2}, Theorem~3.3.)

\begin{lem}
(omit $i=1$ if the corresponding pair consists of integers.)

\medskip

\begin{center}
\begin{tabular}{|lc|} \hline
Transformation & Pairs of resulting branch \\ \hline
 & \\[-6pt]
MONOIDAL TRANSFORMATION & \\ 
$\quad$ Center $(x,z)$ & $\la_i-1,\mu_i$ \\
$\quad$ Center $(y,z)$ & $\la_i,\mu_i-1$ \\
QUADRATIC TRANSFORMATION & \\
``Transversal Case'' ($\la_1+\mu_1\ge 1$) & \\
$\quad$ Direction $(1:0:0)$ & $\la_i+\m_i-1$, $\m_i$ \\
$\quad$ Direction $(0:1:0)$ & $\la_i$, $\la_i+\m_i-1$ \\
``Non-Transversal Case'' ($\la_1+\mu_1<1$) & \\
$\quad$ Direction $(1:0:0)$ &
  $\dis\la_i+\frac{(1+\mu_i)(1-\la_1)}{\mu_1}-2,\frac{1+\mu_i}{\mu_1}-1$ \\
$\quad$ Direction $(0:1:0)$ &
  $\dis\mu_i+\frac{(1+\la_i)(1-\mu_1)}{\la_1}-2,\frac{1+\la_i}{\la_1}-1$ \\
$\quad$ Direction $(0:0:1)$ &
  $\dis\frac{\la_i(1-\mu_1)+\mu_i\la_1}{1-\la_1-\mu_1},
\frac{\la_i\mu_1+\mu_i(1-\la_1)}{1-\la_1-\mu_1}$ \\[3pt] \hline
\end{tabular}
\end{center}

\label{forms}
\end{lem}

In the case of quadratic transformation, we will call the directions
$(1:0:0)$, $(0:1:0)$, and $(0:0:1)$ the $x$-, $y$-, and $z$-charts
respectively.


\section{Embedded Resolution} \label{main}

Let $f(x,y,z)$ define a quasi-ordinary singularity germ at the origin,
and let $\{(\la_i,\mu_i)\}$ be
the collection of characteristic pairs for $f$. We describe the embedded
resolution of f which results from applying the canonical resolution
algorithm of \cite{BM1} to $f$. We show that this process depends solely
on the characteristic pairs of $f$. We also describe the process explicitly;
from our description the characteristic data could be fed to a computer
to output a resolution graph. To describe the process we divide the problem
into two natural stages. First we show how to reduce the multiplicity of
$f$ for the first time; we call this ``the base case.'' Then we show how
to reduce the multiplicity starting from a ``general configuration.''
For this we define a configuration of exceptional divisors arising
after the multiplicity has decreased, and show that it is stable.

\subsection{Outline of Bierstone-Milman's canonical resolution}
The centers for blowing-up an arbitrary germ $f$ are chosen by a process
that can be briefly outlined as follows. Using the Weirstrass Preparation
Theorem, $f$ is written as a pseudo-polynomial of degree $m$ in a
variable $z$, where $m$ is the
multiplicity of $f$. It can also be arranged that the coefficient of the
$z^{m-1}$ term vanishes. Then the surface $N_1=\{z=0\}$
defines a ``maximal contact-space'' for $f$.  The coefficients of $f$
are regarded as functions on $N_1$, one dimension lower than the ambient
dimension. From this collection, a (weighted) ``higher order multiplicity''
is defined and a function is selected which achieves this multiplicity,
which sets up an induction on decreasing dimension. As the
procedure is iterated, the accumulated exceptional divisors are factored out
or included in the coefficient set, depending on when they occurred in the
history of the resolution. Because the codimension cannot increase
indefinitely, the process must end. It does so in one of two ways:
\begin{enumerate}
\item At some stage, there are no more coefficient functions. Then the last
maximal contact space is contained in the vanishing locus of the previous
set of coefficient functions. The final multiplicity is infinity, the order of
vanishing of 0. Since the contact space is automatically smooth, it becomes
the natural choice of center for the next blow-up.

\item A multiplicity of zero is obtained. This occurs when the set of
coefficient functions is generated by a product of exceptional divisors
(which are factored out before the
multiplicity is computed; zero is the order of vanishing of 1). In the
language of \cite{A}, \cite{H}, the coefficient ideal is
principalized.  Each component of the maximum locus of the invariant is an
intersection of exceptional divisors, which is smooth by the inductive
requirement of normal crossings. There is a rule, based on the order in
history of the exceptional divisors, which determines a component to select
as the next center.
\end{enumerate}

The details of this procedure are organized into an invariant $\inv_f$ which
consists of a finite sequence of non-negative rational numbers. The terms of
the sequence are just the weighted multiplicities $\nu_i$ described above
and the numbers $s_i$ of exceptional planes occurring in stages of history
between which the invariant decreased. The numbers $s_i$ reflect the
``age'' of exceptional planes; for example, $s_1$ counts the oldest planes,
those which occurred before the first multiplicity $\nu_1$ decreased.

Ordered lexicographically, $\inv_f$ is an upper-semi-continuous function on
the points of $\{f=0\}$. Its maximum locus defines smooth components which
can be ordered to define canonical centers for blowing up. All of these facts
are explained and proved in \cite{BM1}, \cite{BM2}.

For quasi-ordinary singularities, we define a ``general
configuration of exceptional divisors'', and show that it is stable under every
transformation that is needed in the resolution. (The transformations we
use are blow-ups of points and curves, permuting coordinates, and collecting
integral powers, \ie the transformations used in Lemmas~\ref{trans} and
\ref{forms}). We also show that whenever the multiplicity decreases and
the strict transform is normalized, the resulting collection of exceptional
divisors constitutes a general configuration.

The resolution will be presented in two stages; the first stage is the base
where $E$, the collection of exceptional hyperplanes, is empty. The second
stage is the case where the multiplicity has just dropped, and $E$ is the
general configuration. The first stage sets up the second; we will also use
calculations in the first stage to introduce the machinery of
Bierstone-Milman's canonical embedded resolution.

\bigskip

\subsection{Stage One}

Let $f=z^m+\dotsb +c_0(x,y)$ be a quasi-ordinary germ in
normalized form, and let $(\la,\mu)$ be the first characteristic pair. Assume
this is the beginning of the resolution and $E=\emptyset$. By a series of
canonically-chosen centers, we blow-up $f$ until its multiplicity has been
lowered.  We distinguish two cases: (i) $\la+\mu<1$ (the ``non-transverse
case'' according to Lipman) and (ii) $\la+\mu\ge 1$ (the ``transverse
case''). It is easy to see that the multiplicity $\nu_1$ of $f$ is
$m(\la+\mu)$ in case (i) and $m$ in case (ii).

\bigskip

\noindent {\bf Case (i)} $\quad \la+\mu<1$

Since $f$ is normalized, $\la\mu\ne0$. It follows immediately that the
multiplicity of $f$ at the origin is strictly greater than at any nearby
point. Since $\nu_1$  is the first entry in $\inv_f$, the invariant is
maximized at the origin.  Thus in this case, the canonical center
is immediately determined by the first characteristic pair: the center is the
origin.  In particular it's not necessary to compute the entire sequence
$\inv_f$ to determine the center in this case.

We blow up the origin and examine three coordinate charts where the strict
transforms are determined by $x^{-\nu_1}f(x,xy,xz)$, $y^{-\nu_1}f(xy,y,yz)$,
and $z^{-\nu_1}f(xz,yz,z)$, respectively. We will call them $x$-, $y$- and
$z$-charts. (Lipman used the term $(1:0:0)$, $(0:1:0)$ and $(0:0:1)$
directions, respectively in \cite{L2}.)

In the $x$-chart we have
\[ \begin{array}{rcl}
x^{-\nu_1}f(x,xy,xz) & = & x^{-\nu_1}(x^mz^m+\dotsb +c_0(\zeta_i(x,xy)) \\
 & = & x^{m-\nu_1}z^m+ \dotsb + y^{\nu_1 m}\cdot \mbox{unit}.
\end{array} \]
Using Lipman's formula (Lemma \ref{forms}), the multiplicity has been
lowered to $\nu_{1}'=\mu m$.
To restore $f$ to (perhaps non-normalized) quasi-ordinary form, exchange
variables $y$ and $z$ and apply the preparation theorem.  We obtain
\[ f' = z^{\mu m} + \dotsb + {c'}_0({\zeta'}_i), \]
where the new data is:
\[ (\la',\mu')=\left(\frac{1-\la-\mu}{\mu},\frac{1}{\mu}\right) \]
\[ {\nu'}_1 = \mu m<(\la+\mu)m=\nu_1. \]
Note also that $\mu'>1$ and so $f'$ is in the transverse case. $f'$ is again
a pseudo polynomial defining a quasi-ordinary singularity, but it may not be
normalized: a branch may contain integer powers. To normalize $f'$, replace
$z$ by $z-p(x,y)$ where $p$ is the (convergent) power series of all terms in
$\zeta'$ with integral exponents. Assume for now that $\la'$, $\mu'$ are
not both integral, \ie that $f'$ does not lose a characteristic pair.
It follows from the characterization of
quasi-ordinary parametrization that $x^{\la'}y^{\mu'}|p(x,y)$. The final
configuration is
\begin{equation} f'=z^{m'}+\dotsb+(x^{\la'}y^{\mu'})^m\cdot u,\quad E=\{ (y)\}
\label{config}
\end{equation}
and the multiplicity has dropped.

The $y$-chart is exactly parallel. In the $z$-chart, several blow-ups of the
origin are required (their number depends only on $(\la,\mu)$).
Using Lemma \ref{forms} again, $m$ and $(\la,\mu)$ transform to $m'=m-\nu_1$
and $(\la',\mu')=(\frac{\la}{1-\la-\mu},\frac{\mu}{1-\la-\mu})$; so long as
the new $(\la',\mu')$ satisfies $\la'+\mu'<1$, the multiplicity remains
unchanged. The strict transform is automatically in quasi-ordinary form
and the multiplicity is still maximized at the origin, which
must be the next center. Iterating this process, we continue to blow up
origins in the new $z$-charts until $\la'+\mu'\ge 1$. This last condition
occurs after $\displaystyle k=\left[\frac{1}{\la+\mu}\right]$ blow-ups, where
$[x]$ is the greatest integer function. The final configuration is
\[ f'=z^{m'}+\dotsb +x^{\la m}y^{\mu m}\cdot\mbox{unit},\quad E=\{(z)\} \]
and $m'=m(1-k(\la+\mu))={\nu'}_1$. By the definition of $k$, ${\nu'}_1<\nu_1$.
Using Lipman's formula,
\[ (\la',\mu')=\left(\frac{\la m}{m'},\frac{\mu m}{m'}\right)=
   \left(\frac{\la}{1-k(\la+\mu)},\frac{\mu}{1-k(\la+\mu)}\right).  \]
So $\la'+\mu'>1$, and $f'$ is in the transverse case.

In $x$- or $y$-charts occurring after the first $z$-chart, the multiplicity
drops after just one blow-up. The final configuration is the same as
(\ref{config}) except $E$ also contains $(z)$.

\bigskip

\noindent {\bf Case (ii)} $\quad \la+\mu\ge 1$

We have in this case $\nu_1=m$. We assume that $\la\mu\ne 0$. (The case where
$\la\mu=0$ is simpler.) We assume that $f$ is normalized and $\zeta$
contains no terms with both its exponents integral.
It follows that the coefficient of $z^{m-1}$ in $f$ is zero.
Therefore in this case, a normalized quasi-ordinary germ has a smooth maximal
contact space given by the hyperplane defined by $z=0$.
We will illustrate the calculation of $\inv_f$ and the determination of centers
for a few steps, and then summarize the resolution for this case of stage one.

Proceeding inductively on dimension, let ${\cal H}_1$ be the collection of
coefficients of $f$, regarded as functions on $N_1=\{z=0\}$. More precisely,
${\cal H}_1$ is the set of pairs $(c_q,m-q)$ where $c_q$ is the coefficient
of $z^q$.
The number $m-q$ is the weighted multiplicity assigned to $c_q$. Following
\cite{BM1}, define
\[ \mu_2:=\min\left\{\frac{\mu(h)}{\mu_h}\mid
       (h,\mu_h)\in {\cal H}_1\right\}  \]
where $\mu(h)$ is the order of $h$ at $O\in N_1$. While $E=\emptyset$ we define
$\nu_2=\mu_2$; $\nu_2$ is the ``second multiplicity'' at $O$. Now let
\[ {\cal G}_2=\{(h_1,\nu_2\cdot\mu_h)\mid (h,\mu_h)\in {\cal H}_1\}. \]
Then there must be a pair $(g_{*},\mu_{g_{*}})\in {\cal G}_2$ such that
$\mu(g_{*})=\mu_{g_{*}}$. From here we
repeat the construction:  we find a (linear) change of coordinates which
leaves $z$ unchanged such that a second variable defines a maximal contact
space $N_2\subset N_1$ for the function $g_{*}$. Restricting the coefficients
of $g_{*}$ to $N_2$, define ${\cal H}_2$ and $\nu_3$.

The salient points of this construction are: (see \cite{BM1} and \cite{BM2}
for proofs and details)
\begin{itemize}
\item The numbers $\nu_i$ are analytic invariants of the singularity.
\item The construction at  latter steps of the resolution is amended to
reflect the presence of exceptional divisors.
\item The collection ${\cal G}_i$ can be replaced by an equivalent collection
to make calculations simpler; an ``equivalent'' collection defines the same
value for the invariant in a way that is stable after blowing-up. If
$(g,\mu_g)\in {\cal G}_2$ satisfies $\mu(g)=\mu_g$, and $g$ splits into factors
$g=\prod g_{i}^{m_i}$, then $(g,\mu_g)$ can be replaced by the pairs
$(g_i,\mu_{g_i})$, $\mu_{g_i}=\mu(g_i)$, to obtain an equivalent collection.
\end{itemize}
Now we complete the calculation of $\inv_f$ in the current setting of
$\la+\mu\ge 1$ and $E=\emptyset$. We have
\[ f=z^m+c_{m-2}(x,y)z^{m-2}+\dotsb +c_0(x,y)=\prod_{i=1}^m(z-\zeta_i). \]
Each $c_q$ is a symmetric function of the $\zeta_i$. We have
\[ c_0(x,y)=\prod_{i=1}^m\zeta_i(x,y)=x^{m\la}y^{m\mu}\cdot \mbox{unit}, \]
and
\[ c_q(x,y)=(x^{\la}y^{\mu})^{m-q}\cdot e_q(x,y), \]
where $e_q$ is a fractional power series, possibly a non-unit. This is a
general phenomenon: for any quasi-ordinary singularity, $c_0$
``divides'' the other coefficients in a weighted sense: $h_1$ ``divides''
$h_2$ if $h_{1}^{\mu_{h_2}}\mid h_{2}^{\mu_{h_1}}$. So, $c_0$ determines
the next weighted multiplicity. Therefore,
\[ \nu_2=\mu_2=\frac{m(\la+\mu)}{m}=\la+\mu. \]
The pair $(g,\mu_g)=(c_0, \nu_2\cdot m)\in {\cal G}_2$ satisfies $\mu(g)=\mu_g$,
and since $c_0$ factors, we can replace it by the pairs $(x,1)$ and $(y,1)$.
(Here we use the assumption $\la\mu\ne0$.)
The function $(y,1)$ already defines a second maximal contact space:
$N_2=(y,z)\subset N_1=(z)$. As before, ${\cal H}_2$ is the collection of
coefficient functions of the elements of ${\cal G}_2$, written as polynomials
in $y$. Thus $(x,1)\in {\cal H}_2$, and all other functions $(h,\mu_h)$
of ${\cal H}_2$ satisfy $\mu(h)\ge \mu_h$. Therefore, $\nu_3=\mu_3=1$.
There are no more variables left to define coefficient functions, and $\nu_4$
is defined to be $\infty$. We define $s_i=0$ for $i=1,2,3$ and
\[ \inv_f =(\nu_1,s_1;\nu_2,s_2;\nu_3,s_3;\nu_4)=(m,0;\la+\mu,0;1,0;\infty). \]
When $\inv_f$ ends in the value $\infty$, the last contact space is chosen as
the center. In this case, $N_3=(x,y,z)$ defines the origin. Note that the
invariant depends only  on the characteristic data of $f$, since $m$
is determined by that data.

The  procedure for reducing the multiplicity falls into a handful of cases;
we'll  describe one case completely and then summarize the others.

\bigskip

\noindent ---{\em Year One}

Blow-up the origin and examine the chart defined by $x$.  The strict transform
$f'$ remains in quasi-ordinary form, assuming the multiplicity has not dropped.
Let $O'$ be the origin in this chart. Set $E(O')=\{H_x\}$ where $H_x=(x)$ is
the exceptional divisor. Define $E^1\subset E(O')$ to be the collection of
exceptional hyperplanes
$\{H'\mid H' \mbox{ is the strict transform of } H\in E(O)\}$ where $O$ is the
point that was blown up. Since $E(O)=\emptyset$, $E^1(O')=\emptyset$. In
general, $E^1(O')$ are those exceptional hyperplanes at $O'$ which are pulled
back from the earliest  stage of resolution in which $f$ had the same
multiplicity as it does now.
$\displaystyle E^{i+1}(O')\subset E\setminus\cup^{i}_{q=1}E^q$ is defined
similarly, replacing the multiplicity by the first part of the invariant,
namely the string $\inv_{i-\frac{1}{2}}=(\nu_1,s_1;\dotsb;\nu_i)$.
The numbers $s_i$ are defined as the cardinalities of the sets $E^i$.

For clarity, we'll work through the machinery for calculating $\inv_f$ in
the current chart; after this we give just the essential calculations.

Let ${\cal G}_1=\{(f',m)\}$ be the strict transform of $f$, which we've
assumed retains the same multiplicity. Let
${\cal F}_1={\cal G}_1\cup(E^1,1)$, where $(E^i,1)=\{(g,1)\mid g\in E^i\}$
(Currently, $E^1=\emptyset$ so ${\cal G}_1={\cal F}_1$). From ${\cal F}_1$
we draw a function with $\mu(g)=\mu_g$, and use it to define $N_1$.
Of course, $(f',m)$ is the only choice; since $f'$ is already expressed as
a normalized quasi-ordinary singularity, we may take $N_1=(z)$. As before,
${\cal H}_1=\{({c_i}',m-i)\}$ where the ${c_i}'$ are the coefficients of
$f'$ as a polynomial in $z$. We see that
$\displaystyle\mu_2=\frac{m(\la'+\mu')}{m}=\la'+\mu'$,
where $(\la',\mu')$ is the first characteristic pair for $f'$,
which can be obtained from Lipman's formula (\ref{forms}).

In general, $\mu_2$ is modified to take exceptional divisors into account.
Let $(y_H)=H\in {\cal E}_1=E\setminus E^1$ and define
$\displaystyle\mu_{2H}=\min_{h\in {\cal H}_1}\frac{\mu_H(h)}{\mu_h}$
where $\mu_H(h)$ denotes the order to which $y_H$ factors out of $h$. Then
$\displaystyle\nu_2:=\mu_2-\sum_{H\in {\cal E}_1}\mu_{2H}$. Define
$\displaystyle D_2:=\prod_{H\in {\cal E}_1}{y_H}^{\mu_{2H}}$. By
construction, $D_2^m$ is the greatest common divisor of the elements of
${\cal H}_1$ that is a monomial in exceptional coordinates $y_H$,
$H\in {\cal E}_1$. For each $h\in {\cal H}_1$, write
$h=D_2^m\cdot g$ and $\mu_g:=\mu_h\cdot\nu_2$. Then ${\cal G}_2$ is defined to
be the collection $\{(g,\mu_g)\}$ for all $(h,\mu_h)\in {\cal H}_1$ as above,
together with the pair $(D_2^m,(1-\nu_2)\cdot m)$ if $\nu_2<1$. This
completes a cycle in the definition of $\inv_f$.

Now calculate explicitly: ${c_0}'=(x^{\la'}y^{\mu'})^m\cdot\mbox{unit}$
is the principalizing element of ${\cal H}_1$, and
${\cal E}_1=E\setminus E^1=\{H_x\}$ where $H_x=(x)$. Clearly $\mu_{H_x}=\la'$,
so $\nu_2=\mu'$ and $D_2=x^{\la'}$. We will assume that $\la\ge1$ and
$\mu\ge1$; using Lipman's formula, $\mu=\mu'$, so $\nu_2\ge1$, and $D_2$ is 
not included in ${\cal G}_2$.  We may exclude functions from ${\cal G}_2$
which are divisible by some other function in ${\cal G}_2$. We have
${\cal G}_2=\{(y^{\mu' m},\mu' m)\}\sim \{(y,1)\}$. The first stage of the
invariant, $\mbox{inv}_{1\frac{1}{2}}$, is $(m,0;\mu')=(m,0;\mu)$.
Before blowing up, $\mbox{inv}_{1\frac{1}{2}}=(m,0;\la+\mu)$,
so $\mbox{inv}_{1\frac{1}{2}}$ has decreased.
Therefore $E^2=\{H_x\}$, and ${\cal E}_2$ is defined to be
${\cal E}_1\setminus E^2=\emptyset$. We have $s_2=1$, and
${\cal F}_2={\cal G}_2\cup (E^2,1)=\{(x,1),(y,1)\}$. It's now clear that we
may take $(y,z)$ and $(x,y,z)$ to successively define maximal contact spaces,
and $\nu_3=1$, $s_3=0$, $\nu_4=\infty$. We have
$\inv_f =(m,0;\mu',1;1,0;\infty)$. The last maximal contact space is the
origin, which must be the next center to blow-up.

\bigskip

\noindent ---{\em Year Two}

Blow up the origin, and again examine the $x$-chart. As before, the strict
transform is a normalized quasi-ordinary polynomial. The assumption $\la\ge1$,
$\mu\ge1$ implies the multiplicity cannot decrease by  blowing up the origin,
so $\nu_1=m$. It follows that $E^1=\emptyset$, and ${\cal F}_1={\cal G}_1$
contains only the strict transform, which we will again denote as $f$.
Then as before $N_1=(z)$, and ${\cal H}_1=\{(c_0,m)\}$ (we ignore all higher
coefficient functions,  which are divisible by $c_0$). Again we have
$\mu_2=\la''+\mu''$, where $(\la'',\mu'')$ is the current first characteristic
pair, and ${\cal E}_1=\{H_x\}$ leads to $\mu_{2H_x}=\la''$ and $\nu_2=\mu''$.
Using Lipman's formula, $\mu''=\mu$ and $\mbox{inv}_{1\frac{1}{2}}=(m,0;\mu)$.
This time $\mbox{inv}_{1\frac{1}{2}}$ does not drop, so $E^2=\emptyset$
(the exceptional divisor from the previous stage is blown-away), and
${\cal E}_2={\cal E}_1\setminus E^2=\{H_x\}$. Then
${\cal F}_2={\cal G}_2=\{(y^{\mu m},\mu m)\}\sim\{(y,1)\}$. The next
maximal contact space $N_2$ is defined by $y$ restricted to $N_1$, \ie
$N_2=(y,z)$. There are no coefficient functions contributed to ${\cal H}_2$
by the
only function $(y,1)$ belonging to ${\cal F}_2$, so $\nu_3=\infty$. We have
$\inv_f =(m,0;\mu,0;\infty)$ at this stage. The next center is 
defined by $N_2$, which is the $x$-axis. Note that by assuming $\mu\ge 1$,
the $x$-axis is equimultiple.

\bigskip

\noindent ---{\em Year Three}

Blow-up the x-axis. In this case there is only one chart to consider, the
$y$-chart. (The strict transform of $f$ is smooth in the other chart.)
The strict transform $f'$ of $f$ is $y^{-m}f(x,y,yz)$. $f'$ is (up to
multiplication by a suitable unit) a normalized quasi-ordinary polynomial,
with first characteristic pair $(\la',\mu')=(\la,\mu-1)$ where $(\la,\mu)$
was the corresponding pair in year two. To see the general pattern, we
continue to assume $\la',\mu'\ge 1$. Therefore the multiplicity $m$ has
remained constant. We see that $E=\{H_x,H_y\}$ and $E^1=\emptyset$.
Clearly ${\cal F}_1={\cal G}_1=\{(f',m)\}$, and we take $N_1=(z)$. The
coefficients ${\cal H}_1$ are generated again by
$((x^{\la'}y^{\mu'})^m,m)$, but now
${\cal E}_2=E\setminus E^1=\{H_x,H_y\}$ implies $D_2=x^{\la'}y^{\mu'}$. Then
$\mu_2=\la'+\mu'$ and $\nu_2=0$. The invariant sequence terminates if
$\nu_i=0$ or $\infty$. Thus $\inv_f =(m,0;0)$.

When $\nu_i=0$, the center for the next blow-up is constructed out of the
exceptional locus.  List the elements of $E$ by the order in which each arose
in history; this can be done globally, not just for the current chart.
(Thus we consider $E$ to contain $n$ elements after $n$ blow-ups, most of which
play no role in the calculation of $\inv_f$ at a given point). Any subset $I$
of $E$ can be ordered via the $n$-tuple $(\delta_1,\delta_2.\dotsc\delta_n)$
where $\delta_i=0$ if the {\em i}th element of $E$ is not in $I$; $n$-tuples
are ordered lexicographically.

Returning to the determination of the center of blowing up, when $\nu_2=0$,
define ${\cal G}_2=\{(D_2,1)\}$ and $S_{\mbox{inv}}=\{p\in N_1\mid \mu_p(D_2)
\ge 1\}$.

In the current setting,
\[ S_{\mbox{inv}}=\{p\in N_1\mid \mu_p(x^{\la'}y^{\mu'})\ge 1\}
  =(H_x\cap N_1)\cup(H_y\cap N_1),  \]
the union of the $x$-axis and $y$-axis. The order of $H_x$ is $(0,1,0)$
and the order of $H_y$ is $(0,0,1)$. Thus $H_x$ precedes $H_y$ and the next
center is $C=S_{\mbox{inv}}\cap H_x=\mbox{ the $y$-axis}$.

The procedure now takes on a stable pattern. So long as $\la$ and $\mu$ both
exceed 1, we have $\inv_f=(m,0;0)$. We are forced to continue blowing up axes,
which lowers either $\la$ or $\mu$ by one each time. The corresponding divisor
is blown away and replaced by a newer one, causing the choice of axes to
alternate each time. Eventually one of $\la$ or $\mu$ is driven smaller than
one, and the requirement $\mu_p(D_2)\ge 1$ automatically selects the other
axis ever after, until both $\la<1$ and $\mu<1$. If at this point we still
have still have $\la+\mu\ge 1$, then the multiplicity still remains unchanged.
But the condition $\mu_p(D_2)=\mu_p(x^{\la}y^{\mu})\ge 1$ now selects the origin
as the center, and in fact the origin is the maximum locus for the multiplicity
$\nu_1$. We are forced to blow up the origin repeatedly until $\la+\mu<1$.
(Only $x$-charts and $y$-charts need be considered.)
The multiplicity drops precisely when this occurs. Finally, it is clear that
all of these decisions are completely determined by the arithmetic of the
original pair $(\la,\mu)$.

Except for which charts ($x$ or $y$) we choose to consider in the first two
years and the assumptions on $(\la,\mu)$,
the preceding completely describes how to reduce the multiplicity starting
from year zero in the transversal case. Other starting assumptions
({\em e.g.}, $\la<1$, $\mu\ge 1$) are entirely similar.

\bigskip

\subsection{Stage Two}
Once the multiplicity drops, the algorithm begins again; all
previous history is forgotten. All exceptional divisors passing through the
current point of interest are placed in $E^1$ (the set of ``oldest'' divisors).
When the multiplicity drops for the first time, all the divisors are given by
coordinate functions. In later stages this may not be true. The general
phenomena occurring in later stages is completely illustrated by three
examples.

\begin{eg}
\label{eg1}
Suppose $f$ is normalized, $\la+\mu<1$, and $E=E^1=\{H_x,H_y\}$. This is
the usual situation the first time the multiplicity is lowered when we begin
in the transversal case. Also, note that $\la\mu\ne 0$ since $f$ is normalized.

As before, $\la+\mu<1$ immediately implies the origin has higher multiplicity
than nearby points. We are required to blow up the origin. Examine the
$x$-chart. As before, the multiplicity drops immediately.
The new first characteristic pair is
\[ (\la',\mu')=\left(\frac{1-\la-\mu}{\mu},\frac{1}{\mu}\right), \]
and the strict transform of $f$ has the form
\[ f'=z^mx^{m-\nu}+\dotsb +y^{\mu m}\cdot \mbox{unit}. \]
Note that $\mu'>1$ so we are again in the transversal case. We have
$E=E^1=\{H_x,H_y\}$. To normalize $f'$, we first permute $y$ and $z$,
and, up to a unit factor, get the configuration
\begin{equation}
f'=z^{\mu m}+\dotsb +x^{m-\nu_1}y^m\cdot\mbox{unit}, \quad
  E=\{H_x,H_z\}.
\end{equation}

But now $f'$ may not be normalized; a quasi-ordinary branch $\zeta'$ may
contain terms with integer exponents. Let $p(x,y)$ be the power series
consisting of these terms, and apply transformation (i) of Lemma \ref{trans}.
It follows that $p(x,y)$ is divisible by $x^{\la'}y^{\mu'}$. We now have
$H_{z-p(x,y)}\in E^1=E$. Since $N_1=\{z=0\}$ can be chosen as the
first maximal contact space, both $(x^{\la' m'}y^{\mu' m'},m')$ and
$(p(x,y),1)$ will belong to the coefficient set ${\cal H}_1$, and by
the divisibility relation
between them, the latter function will have no effect on the calculation of
$\inv_f$. This relationship is preserved under blowing-up. Thus the divisor
$H_{z-p(x,y)}$ behaves as though $p(x,y)=0$.
\end{eg}

\begin{eg}
\label{eg2}
Suppose again that $f$ is normalized and $\la+\mu<1$ with
$E=E^1=\{H_x,H_y,H_{z-p(x,y)}\}$ where $x^{\la}y^{\mu}\mid p(x,y)$.
This case may occur when Example \ref{eg1} has lowered it's
multiplicity a second time.

Since $f$ is in the non-transversal case, we must blow up the origin. In the
$z$-chart, the divisor $H_{z-p(x,y)}$ is blown away and plays no further role.
In the $x$-chart, the multiplicity drops immediately. As before, the strict
transform $f'$ requires two steps to be normalized: permute $y$ and $z$,
and then absorb integral powers in the quasi-ordinary branch. In the second
step, we replace $z$ by $z-q(x,y)$ where $x^{\la'}y^{\mu'}\mid q(x,y)$ and
$\displaystyle(\la',\mu')=\left(\frac{1-\la-\mu}{\mu},\frac{1}{\mu}\right)$.
We need to see the effect of these changes on the elements of $E^1$.
Of course $H_x$ is unaffected, and $H_y$ becomes $H_{z-q(x,y)}$, which
behaves as in Example \ref{eg1}. On the other hand, $H_{z-p(x,y)}$ becomes
$H_{y-p(x,z-q(x,y))}$. By the divisibility properties of $p$ and $q$,
$y^{\mu}\mid p(x,y)$ and $x^{\la'}y^{\mu'}\mid q(x,y)$. Since
$N_1=\{z=0\}$ will again be the first contact space, it follows that the
function $y-p(x,z-q(x,y))$ is equal to $y$ up to a unit when restricted
to $N_1$. Thus $H_{y-p(x,z-q)}$ behaves like $H_y$ through out the next
stage of the resolution.
\end{eg}

In Examples \ref{eg1} and \ref{eg2}, we implicitly assumed that Lipman's
formula for $(\la',\mu')$ result in a new first characteristic pair which is
not entirely integral.

\begin{eg}
Suppose $f$ is normalized, $\la+\mu<1$, $E=E^1=\{H_x,H_y,H_{z-p(x,y)}\}$ with
$x^{\la}y^{\mu}\mid p(x,y)$, and $(\la,\mu)=(\frac{a}{m},\frac{1}{m})$ where
$1<a<m$ and $a\in{\mathbb N}$.

As in Example \ref{eg2}, blow up the origin and consider the direction defined
by $x$. The multiplicity of $f'$ drops, but this time
$(\la',\mu')=(m-a-1,m)$. (The new pair becomes integral.) When $f'$ is
normalized, the ``first'' characteristic pair $(\la',\mu')$ has to be
absorbed in the transformation $z\to z-q(x,y)$. Thus
$q(x,y)=x^{\la'}y^{\mu'}\cdot\mbox{unit}$. The strict transform has one fewer
characteristic pairs. (Eventually $f$ has to lose all its characteristic pairs
to become smooth.) 

Just as before, $H_{z}$ becomes $H_{y-p(x,z-q(x,y))}$ and the latter is
generated by $(y)$ once we restrict to $N_1=\{z=0\}$. But while the
transform $H_y\to H_{z-q(x,y)}$ as before, we now have
$f'=z^{m'}+\dotsb +(x^{{\la_2}'}y^{{\mu_2}'})^{m'}\cdot\mbox{unit}$. The
degree-zero coefficient of $f'$,
${c_0}'=(x^{{\la_2}'}y^{{\mu_2}'})^{m'}\cdot\mbox{unit}$, no longer divides
$q(x,y)$. In fact the weighted ``divisibility'' relation is reversed:
$q(x,y)$ ``divides'' ${c_0}'$.
When we construct ${\cal H}_1$, it is $q$ that principalizes the functions and
determines the invariant, until $H_{z-q}$ is blown away. Since the role of $q$
in the process is determined by $(\la,\mu)$, although $(\la',\mu')$ is
integral and no longer counts as a characteristic pair, the process is still
determined by the original characteristic data.  We sketch the next
phase of the resolution. The invariant for the configuration $(f',E)$ is
$(m',3;1,0;1,0;\infty)$ and the next center is the origin. Blow-up the origin
and consider the $x$ direction; the multiplicity can't drop but $H_x$ is blown
away and replaced by a new divisor. The invariant records this as
$(m',2;1,1;1,0;\infty)$, and the origin is the center again. Blow up the
origin and consider the
$y$-direction; now $H_x$ and $H_y$ are ``new'' and belong to ${\cal E}_1$.
The calculation now gives ${\cal H}_1=\{(c_0,m'),(q,1)\}$ and
$D_2=x^{\la'}y^{\mu'}$, where $q(x,y)=x^{\la'}y^{\mu'}\cdot\mbox{unit}$.
We have $\nu_2=0$, and $\inv_f=(m',1;0)$. $D_2$ now determines the centers
for blowing-up by $S=\{p\mid\mu_p(D_2)\ge 1\}$ and the ordering of the
divisors in ${\cal E}_1=\{H_x,H_y\}$. So we are required to blow up the
$x$- and $y$-axis alternately until $\la'<1$ and $\mu'<1$; since
$(\la',\mu')$ are integers, we continue until $\la'=\mu'=0$ and
$H_{z-q(x,y)}$ is blown away.  The relation $q(x,y)\mid c_0(x,y)$ is maintained
though-out, so when $H_{z-q(x,y)}$ is blown away we have
${\cal H}_1=\{(c_0,m')\}$. Then since
$c_0(x,y)=x^{\la_2}y^{\mu_2}\cdot\mbox{unit}$, the second characteristic pair
$(\la_2,\mu_2)$ takes over the process.
\end{eg}

\begin{dfn}
{\em General Configuration}. Let $f$ be a normalized quasi-ordinary
singularity, and let $(\la,\mu)$ be the first characteristic pair.
Furthermore, suppose $f$ was obtained by applying the Bierstone-Milman
algorithm to $\tilde{f}$ with characteristic pairs
$\Lambda=\{(\tilde{\la}_i,\tilde{\mu}_i)\}$. Let
\[ E=E^1\subset\{(z-q(x,y)),(y-q(x,z-p(x,y))),(x)\} \]
 where
\[ y-q(x,z-p(x,y))\mid {}_z=y\cdot\mbox{unit}, \]
and $x^{\la}y^{\mu}\mid q(x,y)$, or $q(x,y)=x^ay^b\cdot\mbox{unit}$ and
$x^ay^b\mid x^{\la}y^{\mu}$, and where $(a,b)$ depends on $\Lambda$.
Then $E$ is called a general configuration.
\end{dfn}

\begin{thm}
At any stage in the resolution when the multiplicity of $f(x,y,z)$ has dropped,
and $f$ is expressed as a normalized quasi-ordinary singularity, the
exceptional divisors constitute a general configuration.
\end{thm}

\pf The theorem is trivial in the base case. We have shown it to be true
at the end of the base case, and in representative examples of the general
case. A complete proof involves checking the calculations for each possible
subset of the general configuration, and seeing that in each case the
final set of divisors (once $f$ is normalized) is again a general
configuration. The calculations are similar in all cases to the three
examples we presented above.
\endpf

Since the controlling function in every case depends only on characteristic
data, we have

\begin{thm}
The Bierstone-Milman canonical resolution algorithm for a quasi-ordinary
singularity $(V,p)$ depends only on the (normalized) characteristic pairs.
\end{thm} 

\begin{eg}
We will find locally an embedded resolution for $(V,p)$ defined by
$f(x,y,z)=z^3+x^2y^4$. The function $f$ has one characteristic pair
$(2/3,4/3)$.

\medskip

\noindent ---{\em Year 0}

In the beginning, $\inf_f(O)=(3,0;2,0;1,0;\infty)$ and the center
is the origin $O$. Blow up $O$ and consider the $x$-chart.

\medskip

\noindent ---{\em Year 1}

$f=z^3+x^3y^4$, $\inf_f(O)=(3,0;4/3,1;1,0;\infty)$ and $E=\{H_x\}$; the center
is again the origin. Blow up the origin and consider the $y$-chart.

\medskip

\noindent ---{\em Year 2}

$f=z^3+x^3y^4$ and $\inf_f(O)=(3,0;0)$ The centers are now
determined by $D_2=xy^{4/3}$ and the ordering of the divisors
$E=\{H_x,H_y\}$. Since $H_x$ is older, we first blow up $H_x\cap N_1$ which
is the $y$-axis.

\medskip

\noindent ---{\em Year 3}

$f=z^3+y^4$, $\inf_f(O)=(3,0;0)$, $D_2=y^{4/3}$. Blow up the
$x$-axis, and consider the $y$-chart.

\medskip

\noindent ---{\em Year 4}

$f=z^3+y$. The multiplicity has dropped, and the strict transform
is now smooth, but the divisors do not have normal crossing with $f$. If
we ``normalize'' $f$ as a quasi-ordinary singularity, we see that $f$ has
lost its only characteristic pair: the previous characteristic pair
becomes $(0,1/3)$, which becomes integral when we normalize by the
transformation $y\leftrightarrow z$, $z\to z+y^3$. We have $f=z$, and
$E=E^1=\{H_x,H_{z-y^3}\}$. Since the multiplicity just dropped, we begin
again, regarding the current state as year 0. We have
$\inf_f(O)=(1,2;1,0;3,0;\infty)$ and the center to blow up is the origin
$O$. Blow up and consider the $x$-chart.

\medskip

\noindent ---{\em Year 1}

In the $x$-chart we have $f=z$ and $E=\{H_x,H_{z-x^2y^3}\}$.
Since of course the multiplicity cannot drop, the new divisor $H_x$ does not
belong to $E^1$. We have $\inf_f(O)=(1,1;3,1;1,0;\infty)$ and the next
center is the origin. Blow up and consider the $y$-chart.

\medskip

\noindent ---{\em Year 2}

$f=z$ and $E=\{H_x,H_y,H_{z-x^2y^4}\}$. Now we have
$E^1=\{H_{z-x^2y^4}\}$ and ${\cal E}_1=\{H_x,H_y\}$, which leads to
$\inf_f(O)=(1,1;0)$ and $D_2=x^2y^4$. All remaining centers are determined
by $D_2$, the exponents of which came from the integral pair $(0,3)$
and the transformations of Lemma~\ref{trans} occurring in the last two
steps. As in the last cycle, we blow up the $y$-axis.

\medskip

\noindent ---{\em Year 3}

$f=z$ and $E^1=\{H_{z-xy^4}\}$, and $\inv_f(O)=(1,1;0)$. Blow up
the $x$-axis and consider the $y$-chart.

\medskip

\noindent ---{\em Year 4}

$f=z$ and $E^1=\{H_{z-xy^3}\}$. Blow up the $y$-axis and consider
the $x$-chart.

\medskip

\noindent ---{\em Year 5--8}

In year 5, $f=z$ and $E^1=\{H_{z-y^3}\}$. Now $D_2=y^3$ and only
the $x$-axis is selected. Blow up the $x$-axis three more times. In year 8,
$f=z$ and $E=\{H_x,H_y,H_{z-1}\}$. We have $\inf_f(O)=(1,0;0)$. The
invariant $\inf_f$ is now locally constant everywhere on $f$, and we have
achieved an embedded resolution.
\end{eg}


\nocite{*}
\bibliographystyle{amsplain}
\bibliography{amsl-bib}

\begin{thebibliography}{DDD}

\bibitem[A]{A}
S.\ S.\ Abhyankar, {\em Resolution of Singularities of Embedded Algebraic
Surfaces}, Springer Monographs in Mathematics (2nd Edition), 1998.

\bibitem[BM1]{BM1}
E.\ Bierstone and P.\ Milman, {\em A simple constructive proof of canonical
resolution of singularities}, Effective Methods in Algebraic Geometry,
Progress in Math {\bf 94} (1991), 11--30.

\bibitem[BM2]{BM2}
\rule[-2pt]{9mm}{.4pt}, {\em Resolution of singularities}, preprint
(alg-geom/9709028).

\bibitem[EV]{EV}
S.\ Encinas and O.\ Villamayor, {\em Good points and algorithmic resolution of
singularities}, Preprint.

\bibitem[G]{G}
Y.-N.\ Gau, {\em Embedded topological classification of quasi-ordinary
singularities}, Mem. Amer. Math. Soc. {\bf 74} (1988), 109--129.

\bibitem[H]{H}
H.\ Hironaka, {\em Resolution of singularities of an algebraic
variety over a field of characteristic zero I, II}, Annals of Mathematics
{\bf 79} (1964), 109--326.

\bibitem[L1]{L1}
J.\ Lipman, {\em Quasi-ordinary singularities of embedded surfaces},
Ph.D. thesis, Harvard University, 1965.

\bibitem[L2]{L2}
\rule[-2pt]{9mm}{.4pt}, {\em Quasi-ordinary singularities of surfaces
in $\C^3$},
Singularities (Proc. Symp. Pure Math. {\bf 40}), Amer. Math. Soc.
Providence 1983, Part 2, 161--171.

\bibitem[L3]{L3}
\rule[-2pt]{9mm}{.4pt}, {\em Topological invariants of quasi-ordinary
singularities},
Mem. Amer. Math. Soc. {\bf 74} (1988), 1--107.

\bibitem[M]{M}
T.T.\ Moh, {\em Canonical uniformization of hypersurface singularities of
characteristic zero},
Journal of Pure and Applied Algebra (to appear).

\bibitem[O]{O}
U.\ Orbanz, {\em Enbedded resolution of algebraic surfaces after Abhyanka
(characteristic 0)}, Lecture Notes in Math., {\bf1101}, Springer, 1984.

\bibitem[V]{V}
O.\ Villamayor, {\em Constructiveness of Hironaka's resolution},
Ann. Scient. Ec. Norm. Sup. $4^e$ serie. t.22 (1989), 1--32

\end{thebibliography}

\bigskip

\noindent {\footnotesize \sc Ohio State University, Mansfield, Ohio 44906}

\noindent {\footnotesize {\it E-mail addresses}:
{\tt cban\verb+@+math.ohio-state.edu},
  {\tt mcewan\verb+@+math.ohio-state.edu}}

\end{document}